\documentclass[12pt]{article}
\usepackage{cancel}
\usepackage{tabularx}
\usepackage{tikz}
\usepackage{relsize}
mm \textheight=8.9in

\usepackage{amssymb}
\usepackage{amsmath}
\usepackage{amsthm}
\usepackage{color}

\usepackage{tikz}
\usetikzlibrary{patterns,snakes}

\usepackage{tikz-cd}

\usepackage{array,multirow}

\renewcommand{\i}{\mathrm{i}}
\newcommand{\br}[3]{{$#1$}$\lower4pt\hbox{$\tp\atop\raise4pt \hbox{$\scriptscriptstyle{#2}$}$} ${$#3$}}
\newcommand{\tw}[3]{{$#1$}${\,\scriptscriptstyle {#2}}\atop\raise9pt\hbox{$\scriptstyle\tp$} ${$#3$}}
\newcommand{\ttps}[2]{{#1}\raise5pt\hbox{$\lower12pt\hbox{$\scriptstyle\tp$}\atop \lower0pt\hbox{$\tilde\;$}$}\raise4.5pt\hbox{${\scriptstyle{#2}}$}}
\newcommand{\st}[1]{\mathrm{st}}

\newcommand{\rd}[1]{\mbox{${\,\scriptscriptstyle {#1}}\atop\raise5.5pt\hbox{$\bullet$}$}}
\newcommand{\rt}[1]{\otimes_\chi}
\newcommand{\lt}[1]{\mbox{${\,\scriptscriptstyle {#1}}\atop\raise5.5pt\hbox{$\ltimes$}$}}
\newcommand{\btr}{\raise1.2pt\hbox{$\scriptstyle\blacktriangleright$}\hspace{2pt}}
\newcommand{\btl}{\raise1.2pt\hbox{$\scriptstyle\blacktriangleleft$}\hspace{2pt}}

\newcommand{\lcr}{\raise1.0pt \hbox{${\scriptstyle\rightharpoonup}$}}
\newcommand{\rcr}{\raise1.0pt \hbox{${\scriptstyle\leftharpoonup}$}}

\newcommand{\ttp}{{\lower12pt\hbox{$\tp$}\atop \hbox{$\tilde\;$}}}

\newcommand{\id}{\mathrm{id}}

\renewcommand{\i}{\mathrm{i}}

\newcommand{\C}{\mathbb{C}}
\newcommand{\Z}{\mathbb{Z}}

\newcommand{\tp}{\otimes}

\newcommand{\zt}{\zeta}

\newcommand{\Fc}{\mathcal{F}}

\newcommand{\gm}{\gamma}
\newcommand{\dt}{\delta}

\newcommand{\End}{\mathrm{End}}

\newcommand{\Rm}{\mathrm{R}}

\newcommand{\La}{\Lambda}

\newcommand{\g}{\mathfrak{g}}

\renewcommand{\k}{\mathfrak{k}}

\newcommand{\h}{\mathfrak{h}}

\newcommand{\s}{\mathfrak{s}}

\renewcommand{\o}{\mathfrak{o}}

\newcommand{\nn}{\nonumber}

\renewcommand{\l}{\mathfrak{l}}

\newcommand{\al}{\alpha}

\newcommand{\bt}{\beta}

\newcommand{\be}{\begin{eqnarray}}
\newcommand{\ee}{\end{eqnarray}}

\newtheorem{thm}{Theorem}[section]
\newtheorem{propn}[thm]{Proposition}

\newcount\prg

\newcommand{\parag}{\advance\prg by1 {\noindent\bf\thesection.\the\prg\hspace{6pt}}}

\begin{document}
\title{Electrical quantum groups}
\author{D. Algethami${}^{\ddag}$, A. Mudrov${}^{\sharp}$, V. Stukopin${}^\sharp$
\vspace{10pt}\\
\small ${\ddag}$ Department of Mathematics, College of Science,\\
\small
 University of Bisha, P.O. Box 551, Bisha 61922, Saudi Arabia,
\vspace{10pt}\\
\small
${\sharp}$ Moscow Institute of Physics and Technology,\\
\small
9 Institutskiy per., Dolgoprudny, Moscow Region,
141701, Russia,
\vspace{10pt}\\
\small
 e-mail: dalgethami@ub.edu.sa,  mudrov.ai@mipt.ru, stukopin.va@mipt.ru
\vspace{10pt}\\
\date{}
}

\maketitle

\begin{abstract}
We give a construction of   electrical quantum group and its superpartner as coideal subalgebras
in $U_q(\g\l(N))$ and, respectively, in $U^u_q(\g\l(1|N-1))$, where the latter
is an "ungraded" (non-standard) quantum general linear supergroup. We derive them as quantum isotropy
subgroups of triangular twisted K-matrices found in one of our recent papers.
\end{abstract}

{\small \underline{2020 Mathematics Subject Classification: 17B37, 17A70}.
\\

\newpage
\section{Introduction}
Reflection equation (RE) is an established topic in quantum integrable models and quantum groups. It was first introduced in \cite{Cher}
and developed further in \cite{Skl,KulSkl,KSS}. Nowadays its theory  is embedded in the theory of coideal subalgebras,
quantum symmetric pairs and their generalizations
\cite{NS,Let,KolbRE}.

In a recent paper of the first two authors \cite{AM}, a classification of solutions to the so called twisted RE associated with the general linear Lie 
superalgebra is obtained. Amongst block-diagonal semisimple K-matrices conjectured for a special polarization in \cite{AMS}, two lower triangular  families 
have been found. One of them consists of even matrices whose classical limit is invariant with respect to the so called electric Lie subalgebra in $\s\l(N)$. 
Such Lie algebras were discovered by Lam and Pylyavsky in their paper on electrical networks \cite{LP}. They were also found to be in connection with  
integrable systems, cluster algebras and Lusztig's positivity problem \cite{GT}. A generalization of electrical Lie subalgebras for general Kac-Moody Lie algebras
is obtained in \cite{BGG}. 

The even triangular K-matrix found in \cite{AM} falls into the standard theory of quantum symmetric pairs of \cite{Let}: there is a coideal subalgebra in the quantum
general linear group stabilizing  $K$ and preserving the RE. One should bear in mined that the corresponding Satake diagram describes
a variety of Lie subalgebras. One of them is the stabilizer of a point on the symmetric space $SL(N)/SO(N)$. Quantization of such a point is an arbitrary invertible diagonal K-matrix.
The lower triangular branch is fixed by the electrical Lie subalgebra. For odd $N$, this is isomorphic to $\mathfrak{sp}(N-1)$. For even $N$, it is the odd symplectic algebra $\mathfrak{sp}(N-1)$, realized as the stabilizer of a non-zero vector in $\mathfrak{sp}(N)$, and has a one-dimensional centre \cite{BGG}, Theorem~1.7.
 Thus the electric subalgebra falls into the standard quantization scheme of symmetric pairs \cite{Let,KolbRE}.

The other  family of triangular K-matrices solving $\Z_2$-graded RE in not even. They are endomorphisms of  a graded vector space $V$ of super-dimension $N-2$ and are relative
to the Lie superalgebrra $\g\l(1|N-1)$. Such solutions exist for  polarizations of $\g\l(1|N-1)$ with at most two odd neighbouring simple roots. 
In the classical limit, these matrices become even and are stabilized by a Lie superalgebra $\mathfrak{k}$. For even $N$, $\mathfrak{k}\simeq\mathfrak{osp}(1|N-2)$. For odd $N$, the even part of $\mathfrak{k}$ is the odd symplectic algebra $\mathfrak{sp}(N-2)$, and $\mathfrak{k}$ has a one-dimensional even centre \cite{AM}, Corollary~9.5.

Contrary to their even counterparts,  non-even K-matrices are generally not associated with  coideal subalgebras in a  quantum supergroup: there is no recipe of constricting
such algebras that would preserve the RE. Non-even  K-matrices can be regarded as 1-dimensional representations of the RE stripped off with its
grading. This suggests the following way of dealing with them. 

It is known that a $\Z_2$-graded Yang-Baxter equation  is equivalent to a non-graded one. The "ungraded" R-matrix is obtained from graded
by the transition $R_{ij,mn}\mapsto R_{ij,mn}(-1)^{\bar i\bar m+\bar m\bar j+\bar j\bar n}$, where $i\mapsto \bar i$ designates the parity of the basis index $i$
equal to the degree of the basis vector $v_i\in \C^N$. It is also known that non-twisted $\Z_2$-graded RE is equivalent to ungraded, see e.g.  \cite{AMS24}, Proposition 2.1.
Equivalence of the  twisted graded and non-graded versions is not known in general. Nevertheless, that is true for the case under the current study. Moreover, with a special choice of the super-transposition, the correspondence between graded and not graded K-matrices can be chosen identical.

Thus we can view our non-even K-matrix as one solving the non-graded twisted RE, whose coideal subalgebra in an appropriate Hopf algebra is available via the standard Noumi-Sugitani dressing
transformation \cite{NS}. 
The total quantum group is non-standard and is related to the original quantum supergroup by the superization procedure of \cite{MRP}.
It is rather close to the ordinary quantum group with a modified Cartan subalgebra action and nilpotence of the isotropic odd root vectors. 
The comultiplication stays similar to the standard case: it is group-like on the Cartan generators and quasi-primitive on simple root vectors.  
The non-standard quantum groups allow for coideal subalgebras that are similar to  Letzter's. One of such 
subalgebras fixes the ungraded K-matrix and preserves the twisted RE. 
In the case of triangular $K$, this can be viewed as a quantization of the electrical Lie superalgebra described in \cite{AM}, Section~9.

\section{Triangular solutions to twisted Reflection equation}
Let   $V$ denote a $\Z_2$-graded vector space $\C^N$ with homogeneous basis $\{v_i\}_{i=1}^N$ of vectors, each of degree $\bar i=\deg(v_i)\in \{0,1\}$ 
The matrix algebra $\End(V)$ is naturally equipped with $\Z_2$-grading returning $\deg(e_{ij})=\bar i+\bar j \mod 2$ on the matrix units. 
The tensor product $\End(V)\tp \End(V)$ is a natural $\Z_2$-graded algebra with multiplication $(A\tp B)(C\tp D)=(-1)^{\deg(B)\deg(C)} AC\tp BD$ for homogeneous  $A,B,C,D\in \End(V)$.

The tensor
$$
R=\sum_{i,j}q^{(-1)^{\bar i}\dt_{ij}} e_{ii}\tp e_{jj}+(q-q^{-1}\sum_{j<i}(-1)^{\bar j}e_{ij}\tp e_{ji}\in \End(V)\tp \End(V)
$$
is known to satisfy $\Z_2$-graded Yang Baxter equation
$$
R_{12}R_{13}R_{23}=R_{23}R_{13}R_{12}.
$$
in $\End(V)\tp \End(V)\tp \End(V)$, see \cite{Zhang}.
The subscripts indicate an embedding of  $R\in \End(V)\tp \End(V)$ in the graded tensor cube $\End(V)\tp \End(V)\tp \End(V)$, 
as generally accepted in  quantum groups.

Let $t$ denote the matrix super-transposition defined on the unit basis matrices by
$$
e_{ij}^t=(-1)^{\bar j(\bar i+\bar j)}e_{ji}.
$$
It is an anti-automorphism of the graded algebra $\End(V)$.
Twisted reflection equation is the identity 
\be
\label{REtw}
R_{21}K_1 R_{21}^{t_1} K_2=K_2 R_{12}^{t_2}K_1 R_{21},
\ee
on a matrix $K\in \End(V)$. The full classification of its solutions is obtained in \cite{AM}.
In particular, it admits two lower triangular solutions
\be
K=\sum_{i=1}^{N}e_{ii} x_i^2+\sum_{i>j=1}^{N} x_i x_j (1+q^{-1}) e_{ij}, \quad  \bar i=0 \quad \mbox{for all} \quad i,
\label{evenK}
\ee
and, for a fixed $m\in \{1,\ldots, N\}$ and $\bar i= \dt_{i,m}$,
\begin{align}
 K
 =\sum_{i=1}^Nx_i^2e_{ii}
 +(1+q^{-1})
   \sum_{\substack{j<i\\i,j\not =m}}
   x_ix_j e_{ij}
 +\eta\sum_{j<m}x_mx_j e_{mj}
 +\eta\sum_{m<i}x_ix_m e_{im}, 
\label{non-evenK}
\end{align}
where $\eta=\sqrt{q^{-1}-q}$. Here $x_i$ are arbitrary non-zero complex parameters.

The matrix (\ref{evenK}) is even while (\ref{non-evenK}) is not. Even K-matrices give rise to coideal subalgebras in quasi-triangular Hopf algebras associated with a matrix $R$,  \cite{NS}. For standard quantum groups, the regular structure theory of coideal subalgebras was developted in  \cite{Let}. Such a theory carries over to 
standard quantum supergroups \cite{AMS}. Non-even solution to (\ref{REtw}) do not fall in the standard theory and require an extended interpretation. 
It will be given with the help of non-standard quantum groups \cite{MRP}.

\subsection{Ungrading Reflection Equation}
\label{Sec_equiv}
The non-graded version of  RE  has the same form as (\ref{REtw}) where the R-matrix, tensor products and the
transposition $t$ are non-graded. In view of equivalence between graded and non-graded Yang-Baxter equations, one can ask if
that equivalence extends to Reflection equations too. The answer is affirmative for  the non-twisted version of RE, see e.g. \cite{Isa} and \cite{AMS24}.
With regard to the twisted version we are not aware of it in general. Yet it is true of the special case of the general linear R-matrix.  
We give a proof for it in this subsection.
 
First we translate operators from $\End(V)\tp \End(V)$ to those in $\End(V\tp V)$:
for homogeneous $A,B\in \End(V)$, with $\bar y$ denoting the parity of $M$,   the graded tensor action reads
$$
 (A\otimes B)^g(v_i\otimes v_j)
       =(-1)^{\bar B \bar i}Av_i\otimes Bv_j.
$$
This is an (non-graded) algebra isomorphism. Specifically,
$$
 R^g=\sum_{i,j}q^{s_i\delta_{ij}}(e_{ii}\otimes e_{jj})^g
       +\omega(q-q^{-1})\sum_{j<i}s_j(e_{ij}\otimes e_{ji})^g,
$$
where $s_i=(-1)^{\bar i}$.
Define  $R^u$ as 
$$
 R^u_{ij,k l}
  =(-1)^{\bar i\bar k+\bar k\bar j+\bar j \bar l}R_{ij,k l}=(-1)^{\bar i\bar k}R^g_{ij,k l}.
$$
This gives the ungraded matrix $R^u$ displayed in  (\ref{eq:R}).
The graded flip is $P^g(v_i\otimes v_j)=(-1)^{\bar i\bar j}v_j\otimes v_i$,
and $R_{21}^g=P^gR^gP^g$.
Define a graded RE residual as a map $\End(V)\ni K\mapsto \Fc(K)\in \End(V)\tp \End(V)$ by the assignment
$$
\Fc(K)=R_{21} K_{1}R_{21}^{t_1}K_{2}
 -K_{2}R^{t_2}_{12}K_{1}R_{21}.
$$
Write $\mathcal F^g(K)$ for the image of $\mathcal F(K)$ in $\operatorname{End}(V\otimes V)$ under the graded tensor action.

Similarly the ungraded residual $\End(V)\ni K\mapsto \Fc^u(K)\in \End(V\tp V)$ is defined as
$$
\Fc^u(K)=R^u_{21} K_{1}(R^u_{21})^{t_1}K_{2}
 -K_{2}(R^u_{12})^{t_2}K_{1}R_{21}^u.
$$
where the $t$ is already the ordinary matrix transposition. Note that  $(R^u_{21})^{t_1}=(R^u_{12})^{t_2}$ because now $t$ is involutive
and $(R^u)^{t_1t_2}=R_{21}^u$.

Introduce an involutive parity operator $J\in \End(V\tp V)$ by $J(v_i\otimes v_j)=(-1)^{\bar i\bar j}v_i\otimes v_j$.
\begin{propn}\label{prop:comparison}
The graded and ungraded Reflection equations are equivalent:
for every  $K\in \End(V)$
\begin{equation}
 \Fc^u(K)=J\,\Fc^g(K)\,J.
 \label{eq:comparison}
\end{equation}
Thus every solution of the graded RE solves the ungraded one and {\em vice versa}.
\end{propn}
\begin{proof}
All matrices in this proof act on the ordinary vector space $V\otimes V$.
 The images of graded R-matrices in $\End(V\tp V)$ are
\be
  R^g&=&\sum q^{s_i\delta_{ij}}e_{ii}\otimes e_{jj}+(q-q^{-1})\sum_{j<i}(-1)^{\bar i\bar j}e_{ij}\otimes e_{ji},
\nn\\ 
(R_{21}^{t_1})^g&=&\sum q^{s_i\delta_{ij}}e_{ii}\otimes e_{jj}+(q-q^{-1})\sum_{j<i}(-1)^{\bar j} e_{ij}\otimes e_{ij},\quad
\nn\\ 
 (R^{t_2})^g&=&\sum q^{s_i\delta_{ij}}e_{ii}\otimes e_{jj}+(q-q^{-1})\sum_{j<i}(-1)^{\bar i} e_{ij}\otimes e_{ij}.
\nn
\ee
The graded flip $P$ commutes with $J$, hence
$R_{21}^g=P^g R^g P^g$.
Directly from these formulas,
$$
 R^u=J R^g,\quad
 R_{21}^u=J R_{21}^g,\quad
 JR_{21}^g=R_{21}^g J,\quad
 (R_{12}^u)^{t_2}=(R_{21}^{t_1})^gJ=J(R_{12}^{t_2})^g.
$$
The middle equality reflects the fact that $R$ is even.
For any $K\in \End(V)$, its graded leg matrices are
$(K_1)^g=K_1$ and $(K_{2})^g=JK_2J$.
Consequently, with $J^2=\id$, we arrive at
\begin{align*}
 R_{21}^uK_1(R_{12}^u)^{t_2}K_2
 &=J\Bigl( R_{21}^g K_1^g (R_{21}^{t_1})^g(JK_2J)\Bigr)J,\\
 K_2(R_{12}^u)^{t_2} K_1R_{21}^u
 &=J\Bigl((JK_2J)(R_{12}^{t_2})^gK_1 R_{21}^g\Bigr)J.
\end{align*}
Subtracting the lines proves the assertion.
\end{proof}

\section{Quantum electrical algebras}
\subsection{Even K-matrix}
Consider the case $\g=\s\l(N)$ first. There are two Satake diagrams with no arcs in this case:
\begin{equation}
\nn
 \begin{picture}(45,10)
\put(2,3){\circle{3}}
\put(3,3){\line(1,0){5}}
\put(11,0){$\cdots$}
\put(26,3){\line(1,0){5}}
\put(32,3){\circle{3}}
 \end{picture}
\quad\quad
 \begin{picture}(45,10)
\put(2,3){\circle*{3}}
\put(3,3){\line(1,0){5}}
\put(11,0){$\cdots$}
\put(26,3){\line(1,0){5}}
\put(32,3){\circle{3}}
\put(34,3){\line(1,0){7}}
\put(42,3){\circle*{3}}
  \end{picture}
\end{equation}
where the one on the right with alternating coloured roots are of odd rank: it has black roots at the ends.
The diagram on the left admits two types of solutions to RE: one is arbitrary diagonal matrix and the other one is the lower
triangular matrix from (\ref{evenK}). 

The diagonal solution is stabilized by a coideal subalgebra which is a quantization
of $\s\o(N)$, it is a "point" on the quantized symmetric space $SL(N)/SO(N)$.
We are interested in a triangular $K$. In the classical limit $q\to 1$ it turns to the matrix $K=\sum_{i=1}^{N}e_{ii} x_i^2+2\sum_{i>j=1}^{N} x_i x_j  e_{ij}$
that is invariant under the electric Lie algebra $\k\subset \g$ generated
by
$$
X_\bt=e_\bt - a_\bt^2  f_\bt- a_\bt[e_\bt,f_\bt], \quad a_{\bt_i}=\frac{x_{i+1}}{x_i}.
$$
For odd $N$, $\mathfrak{k}\simeq\mathfrak{sp}(N-1)$. For even $N$, $\mathfrak{k}$ is the odd symplectic algebra $\mathfrak{sp}(N-1)$, realized as the stabilizer of a non-zero vector in $\mathfrak{sp}(N)$, and has a one-dimensional centre \cite{BGG},  Theorem~1.7. It is non-reductive for even $N\ge4$; for $N=2$, it is one-dimensional and abelian.
 
Quantization of the electrical Lie algebra is a left coideal subalgebra generated by
$$
X_\bt=e_\bt-a_\bt^2 q^{h_\bt} f_\bt -a_\bt \frac{q^{h_\bt}-1}{q-1}, \quad \forall \bt\in \Pi.
$$
The K-matrix (\ref{evenK}) is $\k$-invariant under  twisted adjoint action:
$$
\rho(x) K=K\rho^t(\gm(x)), \quad \forall x\in U_q(\k),
$$
where $\gm$ is the antipode in $U_q(\g)$.

Now we proceed to the study of the case $\g=\g\l(1|N-1)$.

\section{Non-standard (ungraded) quantum group $U_q^{u}(\g)$}

Let $\La$ denote a lattice spanned over $\Z$ by  $Z=\{\zt_1,\ldots, \zt_N\}$ and set $\h^*=\La \tp_\Z \C$. 
Fix a homomorphism  of Abelian groups $\La\to \Z/2\Z$ that  assigns to every basis weight  $\zt$ its parity $\bar \zt$.
Introduce an inner product 
\[
(\zeta_i,\zeta_j)=(-1)^{\bar i}\delta_{ij},\qquad
\]
on the vector space $\La\tp_\Z \C $
and  define a bicharacter $q^{\langle  -, -\rangle}$ on  $\La$ by the formula
$$
q^{\langle \nu,\mu\rangle}=(-1)^{\bar \nu \bar \mu} q^{(\mu,\nu)}, \quad \mu,\nu\in \La,
$$
assuming  $q\in \C^\times $  not a root of unity.

Let $\Rm$ denote the root system of type $A$ with polarization $\Rm=\Rm^-\sqcup\Rm^+$ and the basis $\Pi\subset \Rm^+$ of simple roots,
$$
\Pi=\{\zt_1-\zt_2, \ldots, \zt_{N-1}-\zt_N\}.
$$
Consider   a commutative associative algebra $U_q^u(\h)$ generated by
invertible $T_{\zt}$ with $\zt\in Z$ and put
\[
T_{\sum_i n_i\zeta_i}:=\prod_iT_{\zeta_i}^{\,n_i} \in U_q^u(\h).
\]
In particular, $T_{\alpha_k}=T_{\zeta_k}T_{\zeta_{k+1}}^{-1}$ for each simple root.
Define an associative algebra $U^{u}_q(\g)$ generated over by $U_q^u(\h)$  by simple root vectors $e_\al, f_\al$, $\al\in \Pi$,
subject to:
$$
T_{ \zt} f_\al T_{\zt}^{-1}=q^{-\langle \zt,\al\rangle}f_\al,  \quad T_{ \zt} e_\al T_{\zt}^{-1}=q^{\langle \zt,\al\rangle}e_\al,
$$
$$
[e_\al,f_\bt]=\dt_{\al,\bt}\frac{T_\al-T_\al^{-1}}{q-q^{-1}}
$$
for all $\mu, \nu\in Z$ and $\al, \bt\in \Pi$, and 
to quantum Serre-like  relations
$$
[[e_\al,e_\bt]_{q^{\langle \al,\bt\rangle}}, e_\bt]_{q^{\langle \al+\bt,\bt\rangle}}=0, \quad [[f_\al,f_\bt]_{q^{\langle \al,\bt\rangle}}, f_\bt]_{q^{\langle \al+\bt,\bt\rangle}}=0, \quad (\al,\bt)=\pm 1,
$$
$$
 [e_\al,e_\bt]_{q^{\langle \al,\bt\rangle}}=0 , \quad [f_\al,f_\bt]_{q^{\langle \al,\bt\rangle}}=0 , \quad (\al,\bt)=0,
$$
if $\al\not =\bt$ and
$$
e_\al^2=0, \quad f_\al^2=0, \quad \mbox{if}\quad (\al,\al)=0.
$$
These relations are formulated in terms of  "modified commutator" $[x,y]_c:=xy-cyx$ defined for any invertible scalar $c\in \C^\times$.
There are also higher Serre relations on the root vectors to comply with a  PBW-like theorem for $U^u_q(\g)$, see 
\cite{KwYu}, Definition 2.1, (2.6). We are not using them in this presentation as they automatically fulfilled in the basic representation.

The comultiplication 
$$
\Delta(e_\al)= q^{h_\al }\otimes  e_\al+e_\al\otimes  1,\quad \Delta(T_\mu)=T_\mu\tp T_\mu  ,\quad \Delta(f_\al)=f_\al\otimes  q^{-h_\al } + 1\otimes  f_\al,
$$
where $q^{h_{\al_k}}=T_{\zt_k}T_{\zt_{k+1}}^{-1}$ for all $k=1,\ldots, N-1$, makes $U^{u}_q(\g)$ a Hopf algebra. 
The counit $\epsilon$ vanishes on the root vectors and returns $1=\epsilon(K_\nu)$ on the Cartan generators. 
The antipode $\gamma$ assigns
$$
\gm(K_\nu)=K_\nu^{-1}, \quad \gm(f_\al)=-f_\al q^{h_\al}, \quad\quad \gm(e_\al)= q^{-h_\al}e_\al
$$
for all $\nu\in Z$ and $\al\in \Pi$.

For the constant parity sequence $\bar i=0$, $i=1,\ldots, N$,  $U^{u}_q(\g)$ is the ordinary standard quantum group $U_q(\g)$ of $\g=\g\l(N)$. Otherwise it is related to 
the  quantum supergroup $U_q(\g)$ of general linear superalgebra $\g$ determined by the fixed $\Z_2$-grading of $V$, via a superization procedure \cite{MRP}.
Note that $U^{u}_q(\g)$ is different from the so called bosonization of $U_q(\g)$  obtained via  a $\Z_2$-extension \cite{Yamane}.
 
The algebra $U^{u}_q(\g)$ admits a representation $\rho\colon U^{u}_q(\g)\to \End(V)$ defined via the assignment
$$
\rho(T_{\zt})=\sum_{j=1}^{N} q^{\langle \zt, \zt_j \rangle}e_{jj}, \quad \rho(e_{\al_k})=e_{k,k+1}, \quad \rho(f_{\al_k})=e_{k+1,k}, 
$$
which is compatible with the R-matrix 
\be
\label{eq:R}
 R^{u}=\sum_{i,j}(-1)^{\bar i\bar j}q^{(-1)^{\bar i}\dt_{ij}} e_{ii}\tp e_{jj}+(q-q^{-1})\sum_{j<i}e_{ij}\tp e_{ji}.
\ee

\subsection{Quantization of  electrical superalgebra}
The twisted reflection equation (\ref{REtw}) corresponding to $\g\simeq \g\l(1|N-1)$ has  
valid super-Satake diagrams
\begin{equation} 
\nn
 \begin{picture}(65,10)
\put(2,3){\circle*{3}}
\put(3,3){\line(1,0){5}}
\put(11,0){$\cdots$}
\put(26,3){\line(1,0){5}}
\put(32,3){\circle{3}}
\put(34,3){\line(1,0){7}}
\put(42,3){\circle*{3}}
\put(43.5,3){\line(1,0){8}}
\put(51.5,1.5){\framebox(3,3)}
  \end{picture}
\quad
 \begin{picture}(115,10)
\put(2,3){\circle*{3}}
\put(3,3){\line(1,0){5}}
\put(11,0){$\cdots$}
\put(26,3){\line(1,0){5}}
\put(32,3){\circle{3}}
\put(34,3){\line(1,0){7}}
\put(42,3){\circle*{3}}
\put(43.5,3){\line(1,0){8}}
\put(51.5,1.5){\framebox(3,3)}
\put(55,3){\line(1,0){8}}
\put(63,1.5){\framebox(3,3)}
\put(66,3){\line(1,0){8}}
\put(75,3){\circle*{3}}

\put(76,3){\line(1,0){5}}
\put(84,0){$\cdots$}
\put(99,3){\line(1,0){5}}
\put(105,3){\circle{3}}
\put(107,3){\line(1,0){7}}
\put(115,3){\circle*{3}}
  \end{picture}
\quad
\quad
 \begin{picture}(60,10)
\put(3,1.5){\framebox(3,3)}
\put(6,3){\line(1,0){8}}
\put(15,3){\circle*{3}}

\put(16,3){\line(1,0){5}}
\put(24,0){$\cdots$}
\put(39,3){\line(1,0){5}}
\put(45,3){\circle{3}}
\put(47,3){\line(1,0){7}}
\put(55,3){\circle*{3}}
  \end{picture}
\end{equation}
with at most two odd roots $\al_{m-1}$ and $\al_m$ depicted with the squares \cite{AMSsatake}.
The black-and-white diagrams are of arbitrary odd numbers, and a diagram with an odd root on the right is isomorphic to a  diagram 
with an odd root on the left. 
It gives rise to a family of Lie super algebras $\k\subset\g$ that fix even solutions
to (\ref{REtw}) of a block-diagonal shape with size of blocks at most two. 
The corresponding universal enveloping algebra $U(\k)$ can be quantized  in the standard way as a left coideal subalgebra in $U_q(\g)$:
$$
\Delta(x)\in U_q(\g)\tp U_q(\k), \quad \forall x\in U_q(\k).
$$
The K-matrix is $U_q(\k)$-invariant,
$$
\rho(x) K=K\rho^t\bigl(\zt(x)\bigr), \quad \forall x\in U_q(\k),
$$
as well as  the twisted Reflection equation.

Contrary to the even K-matrices, (\ref{non-evenK}) is not associated with any Satake diagram.  Its classical limit is stabilized by a Lie superalgebra $\mathfrak{k}$ whose even part is an electrical subalgebra in $\mathfrak{sl}(N-1)$. For even $N$, $\mathfrak{k}\simeq\mathfrak{osp}(1|N-2)$. For odd $N$, the even part is the odd symplectic algebra $\mathfrak{sp}(N-2)$, and $\mathfrak{k}$ has a one-dimensional even centre \cite{AM}, Corollary~9.5.
This subalgebra is not quantizable to a coideal
subalgebra preserving  (\ref{REtw}). We suggest the following way of bypassing the obstruction. 

As was demonstrated in Subsection \ref{Sec_equiv}, the matrix $K$ satisfies the "ungraded" RE 
\be
\label{REtwugr}
 R_{21}^{u}  K_1   (R_{21}^{u})^{t_1}   K_2=  K_2  ( R_{12}^{u})^{t_2}  K_1   R_{21}^u,
\ee
with the ordinary  matrix  transposition $t$.
Now one can wrap $K$ up with L-operators of the Hopf algebra $U^{u}_q(\g)$
and get a required  coideal subalgebra in $U^{u}_q(\g)$ stabilizing $K$, as explained in \cite{NS}. This subalgebra  can also  be  described in 
terms of mixed generators, as in ordinary  quantum symmetric pairs. 

Define a subalgebra $U^{u}_q(\k)\subset U^{u}_q(\g)$ generated by $x_i$, $i=1,\ldots, N-1$, where
\be
X_i =
\left \{
\begin{array}{lll}
e_i-\frac{x_{i+1}^2}{x_i^2} q^{h_{\al_i}}  f_i+\frac{x_{i+1}}{x_i}\frac{q^{h_{\al_i}}-1}{\sqrt{q^{-1}-q}}, & i=m-1, m,
\nn\\[10pt]
e_i-\frac{x_{i+1}^2}{x_i^2} q^{h_{\al_i}}  f_i+\frac{x_{i+1}}{x_i}\frac{q^{h_{\al_i}}-1}{1-q}, & i\not =m-1, m.
\end{array}
\right.
\ee
Note with care that $U_q^u(\k)$ is not a result of  ungrading  because a subalgebra $U_q(\k)$ in the quantum supergroup $U_q(\g)$ is not defined. 
\begin{propn}
\begin{itemize}
\item[i)]
$U^{u}_q(\k)$ is a left coideal subalgebra in $U^{u}_q(\g)$.
\item[ii)]
Equation (\ref{REtwugr}) is $U^{u}_q(\k)$-invariant. 
\item[iii)]
The matrix $K$ satisfies
$$
\rho(x) K=  K\rho^t\bigl(\gm(x)\bigr), \quad \forall x\in U^{u}_q(\k).
$$
\end{itemize}
\end{propn}
\begin{proof}
  Straightforward.
\end{proof}
\noindent
Note that $K$ is a non-trivial $U^{u}_q(\k)$-invariant, because it is not fixed by the entire $U^{u}_q(\g)$.  

One can associate the following  Satake diagrams  with $U_q^u(\k)$:
\begin{equation} 
 \begin{picture}(65,10)
\put(2,3){\circle{3}}
\put(3,3){\line(1,0){5}}
\put(11,0){$\cdots$}
\put(26,3){\line(1,0){5}}
\put(32,3){\circle{3}}
\put(34,3){\line(1,0){7}}
\put(42,3){\circle{3}}
\put(43.5,3){\line(1,0){8}}
\put(51.5,1.5){\framebox(3,3)}
  \end{picture}
\quad
 \begin{picture}(115,10)
\put(2,3){\circle{3}}
\put(3,3){\line(1,0){5}}
\put(11,0){$\cdots$}
\put(26,3){\line(1,0){5}}
\put(32,3){\circle{3}}
\put(34,3){\line(1,0){7}}
\put(42,3){\circle{3}}
\put(43.5,3){\line(1,0){8}}
\put(51.5,1.5){\framebox(3,3)}
\put(55,3){\line(1,0){8}}
\put(63,1.5){\framebox(3,3)}
\put(66,3){\line(1,0){8}}
\put(75,3){\circle{3}}

\put(76,3){\line(1,0){5}}
\put(84,0){$\cdots$}
\put(99,3){\line(1,0){5}}
\put(105,3){\circle{3}}
\put(107,3){\line(1,0){7}}
\put(115,3){\circle{3}}
  \end{picture}
\quad
\quad
 \begin{picture}(60,10)
\put(3,1.5){\framebox(3,3)}
\put(6,3){\line(1,0){8}}
\put(15,3){\circle{3}}

\put(16,3){\line(1,0){5}}
\put(24,0){$\cdots$}
\put(39,3){\line(1,0){5}}
\put(45,3){\circle{3}}
\put(47,3){\line(1,0){7}}
\put(55,3){\circle{3}}
  \end{picture}
  \nn
\end{equation}
These diagrams are  not admissible for Lie super-algebras as they lead to improper subalgebras $\k$ exhausting all of $\g$, \cite{AMSsatake}.
Yet  they are  valid for $U^{u}_q(\g)$.

\vspace{20pt}
 
\noindent
\underline{\large \bf Acknowledgement}

\vspace{10pt}
\noindent
This work is done at the Centre of Pure Mathematics MIPT. It is financially supported by Russian Science Foundation grant 26-11-00115.

\noindent
D. Algethami is  thankful to the Deanship of Graduate Studies and Scientific Research at University of Bisha for supporting this work through the Fast-Track Research Support Program.

\section*{Declarations}
\subsubsection*{Data Availability}
 Data sharing not applicable to this article as no datasets were generated or analysed during the current study.

\subsubsection*{Competing interests}
The authors have no competing interests to declare that are relevant to the content of this article.

\end{document}